\documentclass[a4paper,12pt]{amsart}
\usepackage[dvipdfmx]{graphicx}
\usepackage{amsmath,amssymb,amsfonts,amsthm,color}
\newtheorem{definition}{Definition}[section]
\newtheorem{theorem}[definition]{Theorem}
\newtheorem{lemma}[definition]{Lemma}

\newtheorem{example}[definition]{Example}
\newtheorem{proposition}[definition]{Proposition}

\newcommand{\C}{\mathbb{C}}

\newcommand{\M}{\mathbb{M}}

\usepackage{graphicx}
\usepackage{xcolor}
\usepackage{transparent}
\usepackage{import}

\begin{document}

\title{The operator $(p, q)$-norm of some matrices}

\author{Imam Nugraha Albania}
\address[Imam Nugraha Albania]{Mathematics Study Programme, Universitas Pendidikan Indonesia, Jl. Dr. Setiabudhi No. 229, Bandung, West Java, 40154, Indonesia}
\email{albania@upi.edu}
\author{Masaru Nagisa}
\address[Masaru Nagisa]{Department of Mathematics and Informatics, Faculty of Science, Chiba University, 
Chiba, 263-8522,  Japan: \ Department of Mathematical Sciences, Ritsumeikan University, Kusatsu, Shiga, 525-8577,  Japan}
\email{nagisa@math.s.chiba-u.ac.jp}

\maketitle

\begin{abstract}
We compute the operator $(p,q)$-norm of some $n\times n$ complex matrices, which can be seen as bounded linear operators from
the $n$ dimensional Banach space $\ell^p(n)$ to $\ell^q(n)$.
We have shown that a special matrix $A=\begin{pmatrix} 8 & 1 & 6 \\ 3 & 5 & 7 \\ 4 & 9 & 2 \end{pmatrix}$ which corresponds to a magic square 
has $\|A\|_{p,p} = \max \{\|A\xi\|_p : \xi\in\ell^p(n), \|\xi\|_p=1\} =15$ for any $p\in [1,\infty]$.
In this paper, we extend this result and we compute $\|A\|_{p,q}$ for $1\le q \le p \le \infty$.

\medskip\par\noindent
AMS subject classification: Primary 47A30, Secondary 15B05.

\medskip\par\noindent
Key words: $p$-norm, Riesz-Thorin interpolation, $(p,q)$-norm, magic squared matrix.

\end{abstract}

\section{Introduction}
Let $A$ be an $n\times n$ matrix with complex entries.
We call $A$ magic squared if each component is non-negative and
\[   A \xi_0 = \alpha \xi_0 = {}^t A \xi_0 ,  \]
where $\xi_0 = {}^t(1,1,\ldots , 1) \in \C^n$.
We can regard $A$ as a linear operator from $\ell^p(n)$ to $\ell^p(n)$ ($1\le p \le\infty$).
and its norm is defined by
\[    \|A \|_{p,p} = \max \{ \| A \xi \|_p : \xi \in \ell^p(n), \; \| \xi \|_p=1 \} ,   \]
where $\|\xi\|_p = \| {}^t(x_1,x_2,\ldots, x_n) \|_p = (\sum_{i=1}^n |x_i|^p)^{1/p}$.
The method of numerical computation of operator norms treated in  \cite{Boyd} and \cite{Higham}.
Using special properties of matrices, estimations of operator norm was stated in \cite{Bouthat} and \cite{Sahasranand}.
In \cite{Nagisa}, we have shown that, for any magic squared $A$, 
\[   \|A\|_{p,p} = \alpha  \quad \text{for any } p\in [1,\infty].  \]

In this paper, we extend this result to the following statement:

\begin{theorem}
Let $A\in \M_n(\C)$ be magic squared. For $1\le q\le p \le \infty$,
\[    \|A\|_{p,q} =\max \{ \| A \xi \|_q : \xi \in \ell^p(n), \; \| \xi \|_p=1 \}= \alpha n^{\frac{1}{q}-\frac{1}{p}}.   \]
\end{theorem}

\section{Proof of Theorem}

The following two statements are well-known in the interpolation theory (see \cite{Bergh}).

\begin{proposition}[Riesz-Thorin \cite{Bergh}]
Let $A\in \M_n(\C)$ and $p, p_1, p_2, q, q_1,q_2\in [1,\infty]$ with
\[  (\frac{1}{p}, \frac{1}{q}) = (1-\theta)(\frac{1}{p_1}, \frac{1}{q_1}) + \theta (\frac{1}{p_2}, \frac{1}{q_2}) \quad \text{for some } \theta\in[0,1].   \]
Then it holds
\[    \|A\|_{p,q} \le \|A\|_{p_1,q_1}^{1-\theta} \|A\|_{p_2,q_2}^\theta .  \]     
\end{proposition}
\begin{center}
	\includegraphics{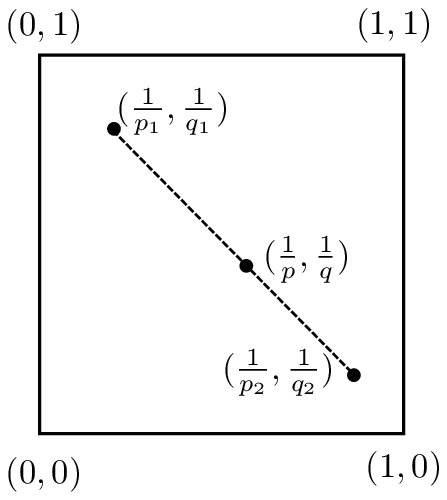}
\end{center}

\begin{lemma}[\cite{Nagisa}:Proposition 3.2]
Let $A\in \M_n(\C)$ and $p_1, p_2, q_1,q_2\in [1,\infty]$.
For any $\theta \in [0,1]$, $p(\theta)$ and $q(\theta)$ is defined as follows:
\[  (\frac{1}{p(\theta)}, \frac{1}{q(\theta)}) = (1-\theta)(\frac{1}{p_1}, \frac{1}{q_1}) + \theta (\frac{1}{p_2}, \frac{1}{q_2}).   \]
If it holds
\[   \|A\|_{p(\theta_0),q(\theta_0)} = \|A\|_{p_1,q_1}^{1-\theta_0} \|A\|_{p_2,q_2}^{\theta_0}  \quad \text{for some } \theta_0\in (0,1) , \]
then
\[   \|A\|_{p(\theta),q(\theta)} = \|A\|_{p_1,q_1}^{1-\theta} \|A\|_{p_2,q_2}^{\theta}  \quad \text{for any } \theta \in [0,1] , \]
\end{lemma}

\vspace{5mm}

\noindent
{\it Proof of Theorem 1.1} \;
At first we prove the following facts:
\begin{enumerate}
  \item[(1)] For any $p\in [1,\infty]$  $\|A\|_{\infty,p} = \alpha n^{1/p}$.
  \item[(2)] For any $p\in [1,\infty]$  $\|A\|_{p',1} = \alpha n^{1/p}$, where $\dfrac{1}{p}+\dfrac{1}{p'}=1$.
  \item[(3)] For any $p\in [1,\infty]$,  $\|A\|_{\frac{2p}{p-1}, \frac{2p}{p+1} }= \alpha n^{1/p}$.
\end{enumerate}
If it holds these 3 properties, we can get
\[   \|A\|_{p(0),q(0)} = \|A\|_{p(1/2),q(1/2)} = \|A\|_{p(1),q(1)} = \alpha n^{1/p},   \]
where
\[   (\frac{1}{p(\theta)}, \frac{1}{q(\theta)}) = (1-\theta)(0, \frac{1}{p}) + \theta (1-\frac{1}{p},1) .  \]
By Lemma 2.2, we have
\[   \|A\|_{p(\theta),q(\theta)} = \alpha n^{1/p} \quad \text{for any } \theta\in [0,1] .   \]
This implies, for any $1\le q \le p \le \infty$, 
\[   \|A\|_{p,q} = \alpha n^{\frac{1}{q}-\frac{1}{p}} , \]
since $(\dfrac{1}{p}, \dfrac{1}{q}) =(1-\theta) (0, \dfrac{1}{q}-\dfrac{1}{p}) + \theta (1-(\dfrac{1}{q}-\dfrac{1}{p}), 1)$ and
$\theta=\dfrac{q}{pq-p+q}$.

So we prove (1), (2), and (3).

(1)  Let $A=(a_{ij})\in \M_n(\C)$ be magic squared and $\xi ={}^t(x_1,x_2,\ldots,x_n) \in \C^n$ with $\|\xi\|_\infty= \max\{|x_i|: i=1,2,\ldots,n\}=1$.
Since $a_{ij}\ge 0$,
\begin{align*}
  \|A\xi\|_p & = \begin{cases}  ( \sum_{i=1}^n | \sum_{j=1}^n a_{ij} x_j |^p )^{1/p}   &  1\le p <\infty  \\
                                        \max \{  | \sum_{j=1}^n a_{ij} x_j | : i=1,2,\ldots, n \}   & p=\infty \end{cases}   \\
             & \le \begin{cases}  ( \sum_{i=1}^n ( \sum_{j=1}^n a_{ij} )^p )^{1/p}   &  1\le p <\infty  \\
                                        \max \{  \sum_{j=1}^n a_{ij}  : i=1,2,\ldots, n \}   & p=\infty \end{cases} \\
             & = \begin{cases}  ( \sum_{i=1}^n  \alpha^p )^{1/p}   &  1\le p <\infty  \\
                                        \max \{  \alpha  : i=1,2,\ldots, n \}   & p=\infty \end{cases} \\
             & = \alpha n^{1/p}.                                                     
\end{align*}
Because $\|A\xi_0\|_p = \alpha n^{1/p}$, we have $\|A\|_{\infty,p} = \alpha n^{1/p}$ for any $p\in [1,\infty]$.

(2) Using the fact that the dual of $\ell^p(n)$ can be identified with $\ell^{p'}(n)$, we have $\|A\|_{p',1} = \|A^*\|_{\infty,p}$, where $\dfrac{1}{p'}=1-\dfrac{1}{p}$.
If we apply (1) for $A^*$ instead of $A$, we can get
\[   \|A\|_{p',1} = \|A^*\|_{\infty,p} = \alpha n^{1/p} .  \]

(3) Since $( \dfrac{1}{\frac{2p}{p-1}}, \dfrac{1}{\frac{2p}{p+1}} ) =\dfrac{1}{2}(0,\dfrac{1}{p}) + \dfrac{1}{2}(\dfrac{1}{\frac{p}{p-1}}, 1 )$, we have
\[    \|A\|_{\frac{2p}{p-1}, \frac{2p}{p+1} } \le \| A\|_{\infty, p} ^{1/2} \| A\|_{p/(p-1), 1}^{1/2} = \| A\|_{\infty, p} ^{1/2} \| A\|_{p', 1}^{1/2} = \alpha n^{1/p}  \]
by (1), (2) and Proposition 2.1.
Since $\|\xi_0\|_{2p/(p-1)} = n^{(p-1)/2p}$, 
\begin{align*}
   \|A\|_{\frac{2p}{p-1}, \frac{2p}{p+1} } & \ge \| A \frac{\xi_0}{ n^{(p-1)/2p}} \|_{2p/(p+1)} = n^{(1-p)/2p} \|\alpha \xi_0\|_{2p/(p+1)}   \\
       & = \alpha    n^{(1-p)/2p} n^{(p+1)/2p} = \alpha n^{1/p}.
\end{align*}
So $\|A\|_{\frac{2p}{p-1}, \frac{2p}{p+1} }= \alpha n^{1/p}$.
\qed

\section{Examples}

For $p\in [1,\infty]$ and $\theta\in (0,1)$, we define $m(p,\theta)$ and $n(p,\theta)$ as follows:
\begin{align*}
   (\frac{1}{m(p,\theta)}, \frac{1}{n(p,\theta)} )& =  (1,0) + \theta \{ ( \frac{1}{p}, \frac{1}{p} ) - (1,0) \}  \\
                  & = ( 1-\theta) (1,0) + \theta (\frac{1}{p}, \frac{1}{p} ) \\
                  & = (1-\theta +\frac{\theta}{p}, \frac{\theta}{p} ) .
\end{align*}
\begin{center}
	\includegraphics{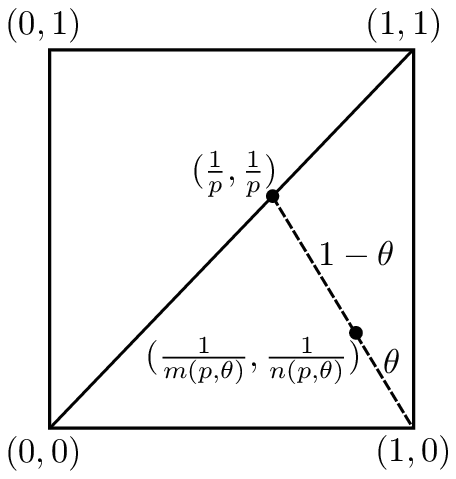}
\end{center}
Then we have
\[   \|A\|_{m(p,\theta), n(p,\theta)} \le \|A\|_{1,\infty}^{1-\theta} \|A\|_{p,p}^\theta  \]
for any $A\in \M_n(\C)$.                 
   
\begin{proposition}
Let $A \in \M_n(\C)$ with $\|A\|_{p,p} = \|A\|_{1,1}$ $(1\le p \le \infty)$.
If $\|A\|_{m(1,\theta_0), n(1,\theta_0)} <\|A\|_{1,\infty}^{1-\theta_0} \|A\|_{1,1}^{\theta_0}$ for some $\theta_0 \in (0,1)$,
then it holds
\[  \|A\|_{m(p,\theta), n(p,\theta)} <\|A\|_{1,\infty}^{1-\theta} \|A\|_{p,p}^\theta  , \]
for any $p\in [1,\infty)$ and $\theta \in (0,1)$.
\end{proposition}
\begin{proof}
At first, we show that the assumption 
\[   \|A\|_{1,\frac{1}{\theta_0}} = \|A\|_{m(1,\theta_0), n(1,\theta_0)} < \|A\|_{1,\infty}^{1-\theta_0} \|A\|_{1,1}^{\theta_0}  \]
implies
\[   \|A\|_{1,\frac{1}{\theta}} = \|A\|_{m(1,\theta), n(1,\theta)} < \|A\|_{1,\infty}^{1-\theta} \|A\|_{1,1}^{\theta}  \]
for any $\theta\in (0,1)$.

We remark that
\[
   (1,\theta) = \begin{cases}  \frac{\theta}{\theta_0}(1,\theta_0) + \frac{\theta_0-\theta}{\theta_0}(1,0)  &  \theta<\theta_0 \\
                     \frac{\theta-1}{\theta_0 -1}(1,\theta_0) + \frac{\theta_0-\theta}{\theta_0 -1}(1,1)  &  \theta>\theta_0 \end{cases}.
\]                
In the case of $\theta < \theta_0$,  we have
\begin{align*}
  \|A\|_{1,\frac{1}{\theta}} & \le \|A\|_{1,\frac{1}{\theta_0}}^{\theta/\theta_0} \|A\|_{1,\infty}^{(\theta_0-\theta)/\theta_0}  \\
       & < (\|A\|_{1,\infty}^{1-\theta_0} \|A\|_{1,1}^{\theta_0})^{\theta/\theta_0}  \|A\|_{1,\infty}^{(\theta_0-\theta)/\theta_0}  \\
       & = \|A\|_{1,\infty}^{1-\theta} \|A\|_{1,1}^\theta .
\end{align*}
In the case of $\theta > \theta_0$,  we also have
\begin{align*}
  \|A\|_{1,\frac{1}{\theta}} & \le \|A\|_{1,\frac{1}{\theta_0}}^{(\theta -1)/(\theta_0 -1)} \|A\|_{1,1}^{(\theta_0-\theta)/(\theta_0 -1)}  \\
       & < (\|A\|_{1,\infty}^{1-\theta_0} \|A\|_{1,1}^{\theta_0})^{(\theta -1)/(\theta_0 -1)} \|A\|_{1,1}^{(\theta_0-\theta)/(\theta_0 -1)}  \\
       & = \|A\|_{1,\infty}^{1-\theta} \|A\|_{1,1}^\theta .
\end{align*}
So $\|A\|_{m(1,\theta), n(1,\theta)} < \|A\|_{1,\infty}^{1-\theta} \|A\|_{1,1}^{\theta} $.

For $p\in [1,\infty)$ and $\theta \in (0,1)$, we assume that
\[   \|A\|_{m(p,\theta), n(p,\theta)} = \|A\|_{1,\infty}^{1-\theta} \|A\|_{p,p}^\theta .  \]
Since $ (\frac{1}{m(p,\theta)}, \frac{1}{n(p,\theta)} ) = (1-\frac{1}{m(p,\theta)})(0,0) + \frac{1}{m(p,\theta)} (1, \frac{m(p,\theta)}{n(p,\theta)} )$, 
we have
\[    \|A\|_{m(p,\theta), n(p,\theta)} \le \|A\|_{\infty,\infty}^{1-\frac{1}{m(p,\theta)}} \|A\|_{1, \frac{n(p,\theta)}{m(p,\theta)}}^{\frac{1}{m(p,\theta)}}.      \]
Since $(1, \frac{m(p,\theta)}{n(p,\theta)}) = \frac{m(p,\theta)}{n(p,\theta)}(1,1) + (1-\frac{m(p,\theta)}{n(p,\theta)}) (1,0)$, we also have
\[   \|A\|_{1,\frac{n(p,\theta)}{m(p,\theta)}} < \|A\|_{1,1}^{\frac{m(p,\theta)}{n(p,\theta)}} \|A\|_{1,\infty}^{1-\frac{m(p,\theta)}{n(p,\theta)}} .  \]
Then this implies the contradiction as follows:
\begin{align*}
     & \|A\|_{1,\infty}^{1-\theta} \|A\|_{p,p}^\theta  =  \|A\|_{m(p,\theta), n(p,\theta)}  \\
     \le  &  \|A\|_{\infty,\infty}^{1-\frac{1}{m(p,\theta)}} \|A\|_{1, \frac{n(p,\theta)}{m(p,\theta)}}^{\frac{1}{m(p,\theta)}}  \\
     <  &   \|A\|_{\infty,\infty}^{1-\frac{1}{m(p,\theta)}}  \|A\|_{1,1}^{\frac{1}{n(p,\theta)}} \|A\|_{1,\infty}^{\frac{1}{m(p,\theta)}-\frac{1}{n(p,\theta)}} \\
     =  &   \|A\|_{p,p}^{1-\frac{1}{m(p,\theta)}+\frac{1}{n(p,\theta)}} \|A\|_{1,\infty}^{\frac{1}{m(p,\theta)}-\frac{1}{n(p,\theta)}}
     = \|A\|_{p,p}^\theta \|A\|_{1,\infty}^{1-\theta},
\end{align*}
where we use the assupmtion $\|A\|_{p,p}=\|A\|_{1,1}=\|A\|_{\infty,\infty}$.
So we have $ \|A\|_{m(p,\theta), n(p,\theta)} < \|A\|_{1,\infty}^{1-\theta} \|A\|_{p,p}^\theta .$ 
\end{proof}

The following fact is well-known as a result of H\"{o}lder's inequality:
\begin{lemma}
Let ${\bf a} = {}^t(a_1, a_2, \ldots, a_n)\in \C^n$ and $p\in [1,\infty]$.
Then we have
\[   \max\{ | \sum_{i=1}^n a_i x_i | : \xi ={}^t(x_1,x_2,\ldots, x_n)\in \ell^p(n), \; \|\xi\|_p =1 \} = \|{\bf a}\|_{p/(p-1)}.  \]
\end{lemma}

\begin{example}
\begin{enumerate}
\item[(1)]  Let $A=I_n= \begin{pmatrix} 1 & & 0 \\ & \ddots & \\ 0 & & 1 \end{pmatrix} \in \M_n(\C)$.
Then
\[   \|A\|_{p,q} = \begin{cases} n^{\frac{1}{q}-\frac{1}{p}}   &  1\le q \le p \le \infty \\  1  & 1\le p<q \le \infty \end{cases}  . \]
\item[(2)]  Let $S=(s_{ij}) \in \M_n(\C)$ and $|s_{ij}| = \begin{cases} 1  & j = \sigma(i)  \\  0 & \text{otherwise} \end{cases}$, where
$\sigma$ is a permutation on $\{1,2,\ldots, n\}$.  
Such a matrix called a unitary permutation.
Then
\[   \|S\|_{p,q} = \begin{cases} n^{\frac{1}{q}-\frac{1}{p}}   &  1\le q \le p \le \infty \\  1  & 1\le p<q \le \infty \end{cases}  . \]
\item[(3)]  Let $A=\begin{pmatrix} 8 & 1 & 6 \\ 3 & 5 & 7 \\ 4 & 9 & 2 \end{pmatrix}$. Then
\[   \|A\|_{p,q} = \begin{cases} 15 \cdot 3^{\frac{1}{q}-\frac{1}{p}},   &   1\le q \le p \le \infty \\
                              (9^{p/(p-1)}+4^{p/(p-1)}+2^{p/(p-1)})^{(p-1)/p}, & 1\le p \le 2, q=\infty \\ 
                              (8^{p/(p-1)}+6^{p/(p-1)}+1)^{(p-1)/p},           & p>2, q=\infty  \\
                              (9^q+5^q+1)^{1/q}                 &  p=1, 1\le q \le \infty    
\end{cases}.   \]
So it holds that $\|A\|_{1,2}=\sqrt{107} <\sqrt{9\cdot 15}=\|A\|_{1,\infty}^{1/2} \|A\|_{1,1}^{1/2}$.

\item[(4)]  Let $a_1,a_2, \ldots, a_n\ge 0$ and $\sigma$ a cyclic permutation on $\{1,2,\ldots,n\}$ of order $n$. 
We define a matrix $A$ as follows:
\[   A = (a_{ij}) \in \M_n(\C), \quad a_{ij} = a_{\sigma^i(j)} \; \text{ for } i,j =1,2,\ldots, n .   \]
Then we have
\[  \|A\|_{p,q} = \begin{cases} (\sum_{i=1}^n a_i) n^{\frac{1}{q}-\frac{1}{p}},  & 1\le q \le p \le \infty  \\
                                       (\sum_{j=1}^n a_j^{\frac{p}{p-1}})^{\frac{p-1}{p}},  &  1\le p \le \infty, q=\infty  \\
                                       (\sum_{j=1}^n a_j^q)^{\frac{1}{q}},   &  p=1, 1\le q \le \infty \end{cases}.  \]
Moreover, in the cases $a=a_1=\cdots=a_n$,
\[    \|A\|_{p,q} = a n^{1-\frac{1}{p}+\frac{1}{q}},  \qquad  1\le p,q \le \infty.  \]                                    
\end{enumerate}
\end{example}

Examples (1), (3), (4) are magic squared. 
So their norm estimations for the case $1\le q\le p \le \infty$ follow from Theorem 1.1.
We explain the rest estimations as below:

(1)  When $p<q$, it is clear that $\| \xi \|_p \le 1$ implies $\| \xi \|_q \le 1$ and
\[   \| \begin{pmatrix} 1\\ 0\\ \vdots \\ 0 \end{pmatrix}\|_p =1 = \| \begin{pmatrix} 1\\ 0\\ \vdots \\ 0 \end{pmatrix}\|_q 
    =\|A\begin{pmatrix} 1\\ 0\\ \vdots \\ 0 \end{pmatrix}\|_q . \] 
So $\|A\|_{p,q}=1$.

(2)  For any $p\in[1,\infty]$ and $\xi\in\ell^p(n)$, $\|S\xi\|_p = \|\xi\|_p$.
In particular, $\|S\|_{p,p}=1$.
Since there exists a unitary permutation $T$ with $TS=I_n$,
\begin{align*}
  \|I_n\|_{p,q} & = \|TS\|_{p,q} \le \|T\|_{q,q}\|S\|_{p,q} = \|S\|_{p,q} \\
     & = \|SI_n\|_{p,q} \le \|S\|_{p,p}\|I_n\|_{p,q} =\|I_n\|_{p,q} . 
\end{align*}
So we have
\[  \|S\|_{p,q} = \|I_n\|_{p,q} = \begin{cases} n^{\frac{1}{q}-\frac{1}{p}}  &  1\le q\le p \le \infty \\ 1 & 1\le p < q \le \infty \end{cases}.  \] 

(3)   For $p\in [1,\infty]$, by using Lemma 3.2,
\begin{align*}
  \|A\|_{p,\infty} & = \max\{ |8x+y+6z|, |3x+5y+7z|, |4x+9y+2z| : \\
                     & \qquad \qquad  |x|^p+|y|^p+|z|^p=1  \} \\
                     & = \max\{ \|{}^t(8,1,6)\|_{p'}, \|{}^t(3,5,7)\|_{p'}, \|{}^t(4,9,2)\|_{p'} \}  \\
                     & =( \max \{ 8^{p/(p-1)}+1+6^{p/(p-1)}, 3^{p/(p-1)}+5^{p/(p-1)}+7^{p/(p-1)}, \\
                     & \qquad \qquad 4^{p/(p-1)}+9^{p/(p-1)}+2^{p/(p-1)} \} )^{(p-1)/p}  \\
                     & = \begin{cases} (9^{p/(p-1)}+4^{p/(p-1)}+2^{p/(p-1)})^{(p-1)/p}, & 1\le p \le 2 \\
                         (8^{p/(p-1)}+6^{p/(p-1)}+1)^{(p-1)/p},  & p>2  \end{cases}.
\end{align*}
In particular, $\|A\|_{\infty, \infty}=15$, $\|A\|_{2,\infty}=\sqrt{101}$ and $\|A\|_{1,\infty}=9$.

Also we have, for $q\in [1,\infty]$ and $q' =\frac{q-1}{q}$, 
\begin{align*}
  \|A\|_{1,q} & = \|A^*\|_{q',\infty}  \\
               & = \max\{ |8x+3y+4z|, |x+5y+9z|, |6x+7y+2z| : \\
               & \qquad \qquad  |x|^{q'}+|y|^{q'}+|z|^{q'}=1  \} \\
               & = ( \max \{ 8^q+3^q+4^q, 1+5^q+9^q, 6^q+7^q+2^q \} )^{1/q}  \\ 
               & = (9^q + 5^q + 1)^{1/q} .    
\end{align*}
In particular, $\|A\|_{1,1}=15$, $\|A\|_{1,2}=\sqrt{107}$ and $\|A\|_{1,\infty}=9$.

(4)  For $p\in [1,\infty]$, by using Lemma 3.2,
\begin{align*}
  \|A\|_{p,\infty} & = \max\{ |\sum_{j=1}^n a_{\sigma^i(j)} x_j | : 1\le i \le n, \sum_{j=1}^n |x_j|^p =1 \}  \\
                     & =(\sum_{j=1}^n a_j^{\frac{p}{p-1}})^{\frac{p-1}{p}} .
\end{align*}

Also we have, for $q\in [1,\infty]$ and $q' =\frac{q-1}{q}$, 
\[  \|A\|_{1. q}  = \|A^*\|_{q', \infty} =(\sum_{j=1}^n a_j^q)^{\frac{1}{q}} . \]

In the case $a=a_1=\cdots=a_n$, we have
\[   \|A\|_{p,q} = a n^{1-\frac{1}{p}+\frac{1}{q}},  \]
since
\[  \|A\xi\|_q \le a\|\xi\|_1 \| \begin{pmatrix} 1 \\ \vdots \\ 1 \end{pmatrix} \|_q \le a n^{1/q} n^{1-1/p}\|\xi\|_p \]
and
\[  \|A \begin{pmatrix} n^{-1/p} \\ \vdots \\ n^{-1/p} \end{pmatrix}\|_q = a n^{1-\frac{1}{p}+\frac{1}{q}}.  \]

\end{document}